\documentclass[10pt,a4paper]{article}
\usepackage{amsfonts}
\usepackage{amssymb}
\usepackage{mathrsfs}
\usepackage{amsthm}
\usepackage{setspace}
\usepackage{color}

\newtheorem{theo}{Theorem}[section]
\newtheorem{prop}[theo]{Proposition}
\newtheorem{cor}[theo]{Corollary}
\newtheorem{lemm}[theo]{Lemma}

\theoremstyle{definition}
\newtheorem*{rem}{Remark}
\newtheorem{defi}[theo]{Definition}

\newtheorem{exa}{Example}

\newenvironment{lproof}{\emph{Proof of Lemma.}}{ \qed \par}

\newcommand{\be}{\begin{eqnarray*}}
\newcommand{\ee}{\end{eqnarray*}}
\newcommand{\beqa}{\begin{eqnarray}}
\newcommand{\eeqa}{\end{eqnarray}}
\newcommand{\ba}{\begin{array}}
\newcommand{\ea}{\end{array}}
\newcommand{\onab}{\overrightarrow{\nabla}}
\newcommand{\mc}{\mathcal}
\newcommand{\mf}{\mathfrak}
\newcommand{\ora}{\overrightarrow}
\newcommand{\rP}{\mathsf{P}}

\newcommand{\mbb}{\mathbb}
\newcommand{\A}{\mc{A}}
\newcommand{\T}{\mc{T}}

\newcommand{\wnab}{\widehat{\nabla}}

\newcommand{\wt}{\widetilde}

\newcommand{\pathgeo}{
\begin{picture}(96,12)\put(5,3){\line(1,0){16}}%
\put(25,3){\line(1,0){16}}%
\put(45,3){\line(1,0){6}}%
\put(71,3){\line(-1,0){6}}%
\put(74,3){\line(1,0){17}}%
\put(59,3){\makebox(0,0){\dots}}%
\put(3,3){\makebox(0,0){$\times$}}%
\put(23,3){\makebox(0,0){$\times$}}%
\put(43,3){\makebox(0,0){$\circ$}}%
\put(73,3){\makebox(0,0){$\circ$}}%
\put(93,3){\makebox(0,0){$\circ$}}%
\end{picture}}

\newcommand{\grasstwogeo}{
\begin{picture}(96,12)\put(5,3){\line(1,0){16}}%
\put(25,3){\line(1,0){16}}%
\put(45,3){\line(1,0){6}}%
\put(71,3){\line(-1,0){6}}%
\put(74,3){\line(1,0){17}}%
\put(59,3){\makebox(0,0){\dots}}%
\put(3,3){\makebox(0,0){$\circ$}}%
\put(23,3){\makebox(0,0){$\times$}}%
\put(43,3){\makebox(0,0){$\circ$}}%
\put(73,3){\makebox(0,0){$\circ$}}%
\put(93,3){\makebox(0,0){$\circ$}}%
\end{picture}}

\newcommand{\skewnil}{\begin{picture}(76,8)\put(5,3){\line(1,0){16}}%
\put(25,3){\line(1,0){6}}%
\put(51,3){\line(-1,0){6}}\put(55,4){\line(1,0){17}}%
\put(54,2){\line(1,0){18}}\put(39,3){\makebox(0,0){\dots}}%
\put(63,2){\makebox(0,0){$>$}}%
\put(3,3){\makebox(0,0){$\circ$}}\put(23,3){\makebox(0,0){$\circ$}}%
\put(53,3){\makebox(0,0){$\circ$}}\put(73,3){\makebox(0,0){$\times$}}%
\end{picture}}

\begin{document}

\title{Generalised Einstein condition and cone construction for parabolic geometries}
\author{Stuart Armstrong, \\ Fakult\"at f\"ur Mathematik, Universit\"at Wien, Nordbergstr. 15, 1090 Wien, Austria}
\date{2008}
\maketitle

\begin{abstract}
This paper attempts to define a generalisation of the standard Einstein condition (in conformal/metric geometry) to any parabolic geometry. To do so, it shows that any preserved involution $\sigma$ of the adjoint Tractor bundle $\mc{A}$ gives rise, given certain algebraic conditions, to a unique preferred affine connection $\nabla$ with covariantly constant rho-tensor $\mathsf{P}$, compatible with the algebraic bracket on $\mc{A}$. These conditions can reasonably be considered the generalisations of the Einstein condition, and recreate the standard Einstein condition in conformal geometry. The existence of such an involution is implied by some simpler structures: preserved metrics when the overall algebra $\mf{g}$ is $\mf{sl}(m,\mbb{F})$, preserved complex structures anti-commuting with the skew-form for $\mf{g}=\mf{sp}(2m,\mbb{F})$, and preserved subbundles of the tangent bundle, of a certain rank, for all the other non-exceptional simple Lie algebras. Examples of Einstein involutions are constructed or referenced for several geometries. The existence of cone constructions for certain Einstein involutions is then demonstrated.

\end{abstract}

\tableofcontents

\section{Introduction}

The study of Einstein manifolds -- spaces whose Ricci tensor is a multiple of the metric -- has been popular ever since Einstein first published his equations \cite{Einstein}, Einstein manifolds being solutions to the vacuum equations of general relativity with a cosmological constant. Many explicit constructions of Einstein spaces have been developed over the years, and it was realised that some metric holonomy groups \cite{berger} force the metric to be Einstein. Conformal and projective parabolic geometries have very close links with Einstein manifolds. When the Ricci-tensor is non-vanishing, these result in an involution of the adjoint Tractor bundle. In this paper, we will find similar constructions for all simple, non-exceptional, parabolic geometries.

Parabolic geometry is a generalisation to non-homogeneous manifolds $M$ of the homogeneous quotient space $G/P$ where $G$ is a semi-simple Lie group and $P$ a parabolic subgroup. The non-homogeneous information is encoded in a principal $P$-bundle $\mc{P} \to M$ and the \emph{Cartan connection}, a one-form $\omega \in \Gamma(\mc{P}, \mf{g})$ for $\mf{g}$ the Lie algebra of $G$. Using a regular Cartan connection to define the geometry, paper \cite{CartEquiv} shows the existence of a unique regular \emph{normal} Cartan connection in that geometry (similar to the way in which the Levi-Civita connection is the unique Torsion-free connection preserving a given metric).

The Cartan connection generates a \emph{Tractor connection} $\onab$ on a principal $G$-bundle $\mc{G}$, which contains $\mc{P}$. Given $\mc{P}$, the Tractor connection $\onab$ and the Cartan connection $\omega$ are equivalent. The standard representation space $W$ of $G$ generates the \emph{standard Tractor bundle}:
\be
\mc{T} = \mc{G} \times_G W,
\ee
on which $\onab$ acts as a vector bundle connection.

Parabolic geometry incorporates many examples of standard geometries. Some of these are given in table \ref{table:17}, characterised by $\mf{g}$ and $\mf{p}$ (it is generally simpler, to avoid issues of connectivity, coverings and centres in $G$, to characterise geometries locally by the Lie algebras rather than globally by the Lie groups). Complexifications and alternative real forms of there geometries are also parabolic geometries, as are many others. The Cartan connection formalism links all of them, but there were few theorems that were truly general. The results of this paper, however, apply to all parabolic geometries, possibly after restricting to an open, dense set of $M$.

\begin{table}[htbp]
\begin{center}
\begin{tabular}{|c||c|c|}
\hline
\hline
Type & algebra $\mathfrak{g}$ & algebra $\mf{p}$ \\
\hline 
\hline
& & \\
Conformal geometry & $\mf{so}(p+1,q+1)$ & $\mf{co}(p,q) \rtimes \mbb{R}^{(p,q)}$ \\
& &\\
Projective geometry & $\mf{sl}(n+1)$ & $\mf{gl}(n)\rtimes \mbb{R}^n$ \\
& & \\
Contact-projective geometry & $\mf{sp}(2n+2, \mbb{R})$ & $( \mbb{R} \oplus \mf{sp}(2n,\mbb{R})) \rtimes \mbb{R}^{2n} \rtimes \mbb{R}$ \\
& & \\
Almost Grassmannian geometry & $\mf{sl}(m + n)$ & $( \mbb{R} \oplus (\mf{sl}(n) \oplus \mf{sl}(m))) \rtimes (\mbb{R}^n \otimes \mbb{R}^m)$ \\
& & \\
CR geometry & $\mf{su}(p+1,q+1)$ & $ (\mbb{R} \oplus \mf{u}(p,q)) \rtimes \mbb{C}^{(p,q)} \rtimes \mbb{R} $ \\
& & \\
Geometry of free $n$-distributions & $\mf{so}(n+1,n)$ & $\mf{gl}(n) \rtimes \mbb{R}^n \rtimes (\wedge^2 \mbb{R}^n)$ \\
& & \\
Path geometry & $\mf{sl}(n+2)$ & $(\mbb{R} \oplus \mf{gl}(n)) \rtimes (\mbb{R} \oplus \mbb{R}^n) \rtimes \mbb{R}^n$ \\
&&\\
\hline
\end{tabular}
\end{center}
\caption{Examples of parabolic geometries}
\label{table:17}
\end{table}

In conformal geometry, $\mf{g} = \mf{so}(p+1,q+1)$, so the standard Tractor bundle $\mc{T}$ has a metric $h$ on it. It was known for a long time \cite{ein} that in this case, a preserved non-degenerate Tractor was locally equivalent to the existence of an Einstein metric in the conformal class. A result by Felipe Leitner \cite{felcon} and by the author \cite{mecon} demonstrated that a preserved, non-degenerate subbundle of $\mc{T}$ implies that the manifold is conformal to a direct product of Einstein manifolds with opposite signs on the Einstein constant. Such subbundles can be characterised by the existence of a preserved metric $g \neq h$ on $\mc{T}$.

Further work by the author in projective geometry \cite{mepro1} and \cite{mepro2} demonstrated that there exists an Einstein connection in the projective class (understood to be an affine connection $\nabla$ that preserves a metric and is Einstein for that metric -- equivalently, that $\mathsf{Ric}^{\nabla}$ is non-degenerate and $\nabla \mathsf{Ric}^{\nabla} = 0$), if and only if $\onab$ preserves a metric $g$ on $\mc{T}$. A similar result was unearthed in the geometry of free $m$-distributions \cite{meskew}: a certain preserved metric $g$ on $\mc{T}$ generates a condition very close to the Einstein condition. This paper aims to generalise this result to any parabolic geometry.

Let $\mc{A}$ be the algebra bundle
\be
\mc{A} = \mc{G} \times_G \mf{g}.
\ee
Via the Cartan connection, there is an inclusion $i: T^* \hookrightarrow \mc{A}$ and a surjective projection $\pi^2: \mc{A} \to T$.

Let $\sigma$ be an involution of $\mc{A}$ -- a map $\sigma: \mc{A} \to \mc{A}$ such that $\sigma^2 = Id_{\mc{A}}$ -- with the algebraic condition that $\pi^2 \circ \sigma \circ i$ is an isomorphism $T^* \to T$. Note that this is automatically the case if $\sigma$ is a Cartan involution (an involution such that $B(\sigma -,-)$ is positive definite).

The main result of this paper is that if there exists such a $\sigma$, this implies the existence of a unique preferred connection $\nabla$ such that the rho tensor $\rP \in \Gamma(\otimes^2 T^*)$ of $\nabla$ is symmetric, non-degenerate, respects the algebraic bracket on $\mc{A}$, and satisfies
\beqa \label{ein:con}
\nabla \rP = 0.
\eeqa
Since $\rP$ is constructed algebraically from the Ricci tensor for normal conformal and projective structures, (and all normal $|1|$-graded geometries) this explains why these structures are referred to as Einstein. The similarity is reinforced by the fact that Equation (\ref{ein:con}) implies that $\nabla$ must be a \emph{metric} connection, using $\rP$ as the metric. The involution $\sigma$ is then called an Einstein involution.

To get these results, there must be ways of dealing with $\mc{T}$ without knowing the details of the subgroup $P$. The main tools to do so is to note that the action of $\mf{g}$ on the standard representation $V$ of $G$ is `nearly transitive' -- specifically, that the span of any non-zero $v \in V$ under the action of $\mf{g}$ is of co-dimension zero, one or two in $V$, depending on $G$. This bundleises to an equivalent statement of the action of $\mc{A}$ on $\mc{T}$, and homogeneity considerations allow analysis of the action of $T \subset \mc{A}$ on $\mc{T}$, without knowing the details of $P$.

These results are very general, but lack one essential ingredient: an existence proof. If we want the Cartan connection to be \emph{normal} (see \cite{TCPG}), the full existence problem can often be simplified. For simple, non-exceptional Lie algebras $\mf{g}$, the existence of such an involution $\sigma$ is implied -- on an open, dense subset of $M$ -- by the existence various simpler structures. A metric on $\mc{T}$ in the cases where $\mf{g}$ is $\mf{sl}(m,\mbb{F})$, a complex structure anti-commuting with the skew form for $\mf{g} = \mf{sp}(2m,\mbb{F})$, and a non-degenerate subbundle of $\mc{T}$ of a certain rank in the case where $\mf{g} = \mf{so}(p,q), \mf{so}(m,\mbb{C}), \mf{so}^*(2m), \mf{su}(p,q),$ or $\mf{sp}(p,q)$. In all cases we ask this extra info be compatible with any complex structures on $\mc{T}$ (commuting or anti-commuting for the metrics and skew-forms).

This implies that, on an open, dense set of $M$, all the holonomy reductions detailed in table \ref{table:one} generate an Einstein involution. Those in table \ref{table:two} imply the existence of an Einstein involution only for certain signatures in the reduced holonomy. Which signatures are valid is dependent on the details of the parabolic inclusion $P \subset G$ -- however for all such geometries, there will be at least one compatible signature generating an Einstein involution. For instance, a preserved subbundle $K \subset \mc{T}$ of rank $r(K) = r(\mc{T})/2$ or $k= r(\mc{T})/2 - 1/2$ \emph{always} generates an Einstein involution in this case.

\begin{table}[htbp]
\begin{center}
\begin{tabular}{|c|c||c|c|}
\hline
\hline
algebra $\mathfrak{g}$ & holonomy reduction & algebra $\mathfrak{g}$ & holonomy reduction \\
\hline 
\hline
&&&\\
$\mf{sl}(m,\mbb{R})$ & $\mf{so}(p,q)$ & $\mf{sp}(2m,\mbb{R})$ & $\mf{u}(p,q)$ \\
&&&\\
$\mf{sl}(m,\mbb{C})$ & $ \mf{su}(p,q), \ \mf{so}(m,\mbb{C})$ & $\mf{sp}(2m,\mbb{C})$ & $ \mf{sp}(p,q), \ \mf{gl}(m,\mbb{C})$ \\
&&&\\
$\mf{sl}(m,\mbb{H})$ & $ \mf{sp}(p,q), \ \mf{so}^*(2m)$ & & \\
&&&\\
\hline
\end{tabular}
\end{center}
\caption{Holonomy reductions implying an Einstein involution, $p+q = m$}
\label{table:one}
\end{table}

\begin{table}[htbp]
\begin{center}
\begin{tabular}{|c|c||c|c|}
\hline
\hline
algebra $\mathfrak{g}$ & holonomy reduction & algebra $\mathfrak{g}$ & holonomy reduction \\
\hline 
\hline
&&&\\
$\mf{so}(p,q)$&$\mf{so}(p',q') \times \mf{so}(p'',q'')$&$\mf{su}(p,q)$&$\mf{su}(p',q')\times\mf{su}(p'',q'')$ \\
&&&\\
$\mf{so}(m,\mbb{C})$&$\mf{so}(p',\mbb{C}) \times \mf{so}(q',\mbb{C})$&$\mf{sp}(p,q)$&$\mf{sp}(p',q')\times\mf{sp}(p'',q'')$ \\
&&&\\
$\mf{so}^*(2m)$&$\mf{so}^*(2p') \times \mf{so}^*(2q')$&& \\
&&&\\
\hline
\end{tabular}
\end{center}
\caption{Holonomy reductions implying an Einstein involution for certain $p',q',p'',$ and $q''$}
\label{table:two}
\end{table}

These various structures become equivalent with the existence of solutions of a series of invariant differential equations. Without attempting to solve these equations explicitly, this paper will instead give existence results for conformal, projective, contact-projective, CR, path, almost quaternionic and almost Grassmannian geometries. Moreover if the Einstein involution is also a Cartan involution, it generates Einstein involutions on all its correspondence spaces (see \cite{CoresSpace} for more details on correspondence space). That can generate many more examples, as, for instance, projective and conformal structures with preserved Cartan involutions exist, and these have a vast amount of correspondence spaces.

The final section deals with a generalisation of the cone construction that exist for projective geometries and conformally Einstein conformal geometries. Though the presence of an Einstein-involution does not guarantee the existence of a cone construction, it does make it more likely that such a construction exists. An example of this construction in the case of conformal, $m$-distribution, path and almost Grassmannian geometries is given in the last section.

\subsection*{Acknowledgements}
It gives me great pleasure to acknowledge the financial support of an ESI Junior Fellowship program and project P19500-N13 of the ``Fonds zur F\"orderung der wissenschaftlichen Forschung (FWF)'', as well as the help, proofreading and comments of Andreas {\v{C}}ap.

\section{Parabolic geometries and metrics}

\subsection{Cartan connections}
This section will present the formalism for Cartan/Tractor connections on parabolic geometries to sufficient depth to set the notions and notations for this paper. See \cite{TCPG} for a good general introduction to parabolic geometries; \cite{CartEquiv} and \cite{two} are also good sources.

A homogeneous space is a space $M = G/P$, where $G$ and $P$ are Lie groups. This makes $G$ into an $P$-bundle over $M$. The left invariant vector fields on $G$ define an isomorphism between the tangent space $TG_g$ for all $g \in G$ and $TG_{id}$, the tangent space at the identity. Since $TG_{id} \cong \mf{g}$, the Lie algebra of $G$, this isomorphism is equivalent with a section $\omega$ of $TG^* \otimes \mf{g}$. It is easy to see that this is $P$-equivariant.

A Cartan connection $\omega$ on a manifold $M$ is a generalisation of this idea to non-homogeneous manifolds $M$. Specifically, it is provided by a principal $P$-bundle $\mc{P} \to M$ and a section $\omega$ of $T\mc{P} \otimes \mf{g}$ with the following properties:
\begin{enumerate}
\item $\omega$ is $P$-equivariant.
\item At any point $u \in \mc{P}$, $\omega_u: T\mc{P}_u \to \mf{g}$ is an isomorphism.
\item If $A \in \mf{p}$, the Lie algebra of $P$, and $\xi_A$ is the vector field on $\mc{P}$ generated by $A$, $\omega(\xi_A) = A$.
\end{enumerate}

The second property shows that $\omega$ is not a connection on a principal bundle in the standard sense (it does not define a horizontal subspace of $T\mc{P}$), and thus cannot be used for differentiating on bundles associated with $\mc{P}$. However, the inclusion $P \subset G$ generates a bundle inclusion $i: \mc{P} \subset \mc{G}$, with $\mc{G} \to M$ a principle $G$-bundle. There is a unique $G$-equivariant section $\omega'$ of $T\mc{G}^* \otimes \mf{g}$ such that $\omega'(\xi_A) = A$ for all $A \in \mf{g}$ and $\omega= i^* \omega'$. This is a connection on the principle bundle $\mc{G}$, the so-called \emph{Tractor} connection.

The Cartan geometry is provided by $\mc{P}$ and $\omega$ (since $\omega$ and $\omega'$ are equivalent given $\mc{P}$, we will suppress the distinction between them). If $V$ is any representation of $G$, we can form the bundle
\be
\mc{V} = \mc{G} \times_G V,
\ee
and $\omega$ generates a connection on $\mc{V}$, designated by $\onab$. Since any representation of $G$ is, a fortiori, a representation of $P$, we have
\be
\mc{V} = \mc{P} \times_P V,
\ee
giving us extra structure on $\mc{V}$. We shall call these bundles -- bundles associated to $\mc{P}$ via the restriction to $P$ of a representation of $G$ -- \emph{Tractor bundles}. The \emph{standard Tractor bundle} is that generated by the standard representation of $G$, and is designated $\mc{T}$. The adjoint Tractor bundle is that generated by the adjoint representation of $G$, and is designated $\mc{A}$.

\subsection{Parabolic geometries} \label{parabolic:section}
A parabolic geometry is one where the inclusion $P \subset G$ is parabolic. There are invariant ways of seeing this property \cite{paradef}, but a simple characterisation will suffice here:
\begin{defi}
A subgroup $P$ of a connected semi-simple Lie group $G$ is parabolic if the Lie algebra $\mf{g}$ of $G$ admits a grading:
\be
\mf{g}_{-k} \oplus \ldots \mf{g}_{-1} \oplus \mf{g}_0 \oplus \mf{g}_1 \oplus \ldots \mf{g}_k,
\ee
such that there are no simple ideals of $G$ in $\mf{g}_0$, $[\mf{g}_i, \mf{g}_j] \subset \mf{g}_{i+j}$ and the Lie algebra of $P$ is 
\be
\mf{p} = \sum_{j\geq 0}^k \mf{g}_j.
\ee
\end{defi}
This grading is not uniquely defined; it changes by the action of $P$. The filtered subspaces $\mf{g}_{(i)} = \sum_{j\geq i}^k \mf{g}_j$ however, are well defined; this will be a general characteristic of structures associated to parabolic geometries. The Tractor bundles, for instance will have a filtration.

Define 
\be
\mc{A}_{(j)} = \mc{P} \times_P \mf{g}_{(j)}.
\ee
These give a filtration of the adjoint Tractor bundle $\mc{A}$ as
\be
\mc{A} = \mc{A}_{(-k)} \supset \ldots \supset \mc{A}_{(0)} \supset \ldots \supset \mc{A}_{(k)}.
\ee
Note that these filtered bundles are not Tractor bundles -- the action of $P$ on $\mf{g}_{(j)}$ does not come from a restriction of the action of $G$.

Since $\omega$ maps the vertical vectors of $\mc{P}$ to elements of $\mf{p}$, we may use $\omega$ to identify the pull back of $TM$ to $\mc{P}$ at each point $u \in \mc{P}$ with $\mf{g}/\mf{p}$. Since $\omega$ is $P$-equivariant, we may divide out by the action of $P$ and get the relation:
\be
\mc{A} / \mc{A}_{(0)} = TM.
\ee
The Killing form on $\mf{g}$ identifies $\mf{g}_j$ with $\mf{g}_{-j}^*$. Passing to the bundle, this implies that
\be
\mc{A}_{(1)} = TM^*.
\ee

We may further define the associated graded bundles $\mc{A}_{j} = \mc{A}_{(j)} /\mc{A}_{(j+1)}$. Gradings are generally easier to handle than filtrations; but the Tractor connection does not operate on these gradings. What we would want is an isomorphism between the graded algebra bundles and the filtered ones. This is done through the choice of a \emph{Weyl structure}:
\begin{defi}
A Weyl structure is given by a filtration preserving algebra isomorphism
\be
\mc{A} = \sum_{j = -k}^k \mc{A}_j,
\ee
\end{defi}
such that $\mc{A}_{(i)} = \sum_{j=i}^k \mc{A}_j$.

There are other ways of looking at Weyl structures, such as the existence of the \emph{grading section}. Since the endomorphism $\theta: \mf{g} \to \mf{g}$, $\xi_j \to j \xi_j$ for $\xi_j \in \mf{g}_j$ is an inner endomorphism, there must exist an element $\wt{e} \in \mf{g}$ such that $ad_{\wt{e}} = \theta$. Since
\be
ad_{\wt{e}} \wt{e} = [\wt{e},\wt{e}] = 0,
\ee
we must have $\wt{e}\in \mf{g}_0$. The above construction implies that $\wt{e}$ defines the grading; it shall be called the grading element.
\begin{lemm}
The image $e$ of $\wt{e}$ under the projection $\mf{g}_{(0)} \to \mf{g}_0$ is the same for all splittings of $\mf{g}$ compatible with $P$.
\end{lemm}
\begin{lproof}
The compatible splittings of $G$ change under the action of $P$, and $\wt{e}$ changes in the same way. Since $ad_{\wt{e}} (\mf{p}) \subset \mf{g}_{(1)}$ acts trivially on $\mf{g}_0$ under the quotient action, the result follows.
\end{lproof}
This unique element $e$ allows us to define a section $E$ of $\mc{A}_0$, $E = \mc{P} \times_P e$, with $P_{(1)}$ acting by trivial quotient action on $e$. Now a Weyl structure is equivalent with a $P$-equivariant map from $\mc{P}$ to $\mf{g}_{(0)}$ that projects to $e$ under the quotient projection $\mf{g}_{(0)} \to \mf{g}_0$. This map defines a section $\wt{E}$ of $\mc{A}_{(0)}$ which we shall call the \emph{grading section}. By construction, it is a lift of $E$ from $\mc{A}_0$ to $\mc{A}_{(0)}$.

If $G_0$ is the subgroup of $G$ whose Lie algebra is $\mf{g}_0$, the projection $P \to G_0$ (dividing out by the group generated by $\mf{g}_{(1)}$) defines a bundle projection $\mc{P} \to \mc{G}_0$. Now, given a Weyl structure, we have a grading of the Tractor connection:
\be
\omega=& \omega_- &(= \omega_{-k} + \ldots + \omega_{-1}) \\
&+ \omega_0 \\
&+ \omega_+ &(=\omega_1 + \ldots + \omega_k).
\ee
Since $\omega_j$ is preserved by the (trivial) quotient action of $\mc{A}_{(1)}$, these $\omega_j$ descend to forms on $\mc{G}_0$. Under our identification $\mc{A} / \mc{A}_{(0)} = TM$, $\omega_-$ is simply the identity on $TM$ (though the particular form of $\omega_- = \omega_{-k} + \ldots + \omega_{-1}$ does give us a grading on $T$). This makes $\mc{G}_0$ into a principal bundle for $TM$. The central term is $G_0$ equivariant and maps vertical elements of $\mc{G}_0$ to $\mf{g}_0$, making it into a principal connection on $\mc{G}_0$, hence an affine connection $\nabla$ on $TM$.

These connections $\nabla$ are called preferred connections, and are equivalent with both the Weyl structure and the grading section $\wt{E}$. Finally the last piece $\omega_+$ is a section of $TM^* \otimes TM^*$ dependent on $\nabla$ and designated $\rP$, the \emph{rho} tensor. Thus we may define, for each preferred connection $\nabla$, a splitting of the adjoint Tractor bundle
\be
\mc{A} = TM \oplus \mc{A}_0 \oplus TM^*
\ee
and express the Tractor connection as:
\be
\onab_X = X + \nabla_X + \rP(X).
\ee
Given a grading section $\wt{E}$, we may split any Tractor bundle $\mc{V}$ into eigenspaces of $\wt{E}$, with eigenvalue $j$. These eigenbundles will be designated $H_j$, and are said to have homogeneity $j$. The action of homogeneous elements of $\mc{A}$ interchanges these bundles. Since the homogeneities of $\mc{A}$ are all integers, if $\mc{V}$ comes from an irreducible representation $V$ of $G$, then the homogeneities of $\mc{V}$ must differ by integers. These bundles do depend on the choice of Weyl structures, but the filtered bundles
\be
\mc{V} \supset H_{(j)} = \sum_{i\geq j} H_i,
\ee
are well defined, independently of $\wt{E}$. Well defined also is the \emph{highest homogeneity subbundle} $H_{(l)} = H_l$.

\section{Preserved involutions: generalised Einstein manifolds}
\subsection{Einstein involutions}

\begin{theo} \label{main:theo}
Let $\sigma$ be an involution of the algebra bundle $\mc{A}$ such that $\onab \sigma = 0$ and
\beqa
\label{intersection:lemm} \sigma(T^*) \cap \mc{A}_{(0)} = 0
\eeqa
(equivalently, $\pi \circ \sigma: T^* \to T$ is bijective). Then the holonomy of $\onab$ is contained in the $+1$ eigenspace of $\sigma$, and there is a unique preferred connection $\nabla$ defined by $\sigma$. This $\nabla$ has the following properties:
\begin{itemize}
\item[-] $\rP$ is non-degenerate and symmetric, hence is a metric on $M$,
\item[-] using the splitting defined by $\nabla$ to decompose $T^* = \sum_{j=1} T^*_j$ into homogeneous components, $\rP$ is a section of $\sum_j T^*_j \otimes T^*_j$.
\item[-] $\{\rP(X), \rP(Y)\} = \rP(\{X,Y\})$, for $X$ and $Y$ sections of $T$,
\item[-] $\sigma$ acts as $\rP: T \to T^*$,
\item[-] $\sigma$ restricts to an involution of $\mathcal{A}_0$, and on that bundle it is minus the action of $\rP$ acting by conjugation on $T \otimes T^*$,
\item[-] $\nabla \rP = 0$.
\end{itemize}
And conversely, any such $\nabla$ defines an involution $\sigma$. In the splitting defined by $\nabla$, the $+1$ eigenspace of $\sigma$ is an algebra bundle generated by elements of the form
\be
X + \rP(X),
\ee
for $X$ any section of $T$. The holonomy algebra of $\onab$ must then reduce to this eigenspace.
\end{theo}
Any involution that obeys property (\ref{intersection:lemm}) is called an \emph{Einstein involution}. Let $F_+$ be the $+1$ eigenspace of $\sigma$, $F_-$ the $-1$ eigenspace. Since $\sigma$ is an involution, $\mc{A} = F_+ \oplus F_-$. Since $\sigma$ preserves the Lie bracket, $F_+$ is an algebra bundle. Moreover, both $F_+$ and $F_-$ are of locally constant rank -- this can be seen by parallel transport using $\onab$, which, since $\onab \sigma = 0$, must preserve $F_+$ and $F_-$.

Proof of this theorem will come from the following two lemmas:
\begin{lemm}
There is an subbundle $\mc{C}$ of $\mc{A}_{(0)}$ that projects bijectively onto $\mc{A}_0$ such that $\sigma(\mc{C}) = \mc{C}$. This allows us to defined an involution $\sigma$ on $\mc{A}_{0}$.
\end{lemm}
\begin{lproof}
Let $a$ be the rank of $\mc{A}_0$, and $n$ the dimension of $M$. By definition, $\mc{A}$ is of rank $a+2n$, $\mc{A}_{(0)}$ of rank $a+n$ and $\mc{A}_{(1)} = T^*$ of rank $n$. Let $r_1$ be the rank of $F_+$, and $r_2 = 2n + a - r_1$ the rank of $F_-$. By equation (\ref{intersection:lemm}), $r_1$ and $r_2$ are less than or equal to $n+a$ (otherwise, they would have an intersection with $T^*$, giving a section of $T^*$ stabilised by $\sigma$). This also implies that they are both greater than or equal to $n$.

Define
\be
\mc{C} = \mc{C}_+ \oplus \mc{C}_-,
\ee
where $\mc{C}_+ = (F_+ \cap \mc{A}_{(0)})$ (of rank $r_1 - n$) and $\mc{C}_- = (F_- \cap \mc{A}_{(0)})$ (of rank $r_2-n$). Consequently $\mc{C}$ is of rank $r_1 + r_2 - 2n = a$. We now need to show that the projection $\pi_0: \mc{A}_{(0)} \to \mc{A}_{0}$ projects $\mc{C}$ bijectively onto $\mc{A}_{0}$.

Let $t \in \Gamma(\mc{C})$ be a local section such that $\pi_0(t) = 0$. This means that $t$ is a section of $T^*$ and $t = t_+ + t_-$, where $t_+$ and $t_-$ are sections of $\mc{C}_+$ and $\mc{C}_-$ respectively. Applying $\sigma$ to $t$ defines $\sigma(t) = t_+ - t_-$, a section of $\mc{A}_{(0)}$. Then equation (\ref{intersection:lemm}) implies that $t = 0$.

Since $\pi_0$ is an algebra homomorphism, $\mc{C} \cong \mc{A}_0$ is an algebra bundle and $\sigma$ descends to an involution of $\mc{A}_0$.
\end{lproof}

Now consider the algebra $\mf{g}_0$ with an involution $s$ on it. Let $\xi$ be any element of $\mf{g}_0$ and $a$ any element of the centre of $\mf{g}_0$. Then
\be
0 = s [\xi, a] = \pm [ s(\xi), a].
\ee
Thus $s$ preserves the centre of $\mf{g}_0$, and, shifting to the bundle point of view, $\sigma$ preserves the centre of $\mc{C}$.

Let $E$ be the grading section of $\mc{A}_0$. We can lift $E$ to the corresponding grading section $\wt{E}$ of $\mc{C} \subset \mc{A}_{(0)}$. This gives us a Weyl structure, hence a preferred connection $\nabla$ and a splitting of $\mc{A}$. Since all of $\mc{C}$ commute with $\wt{E}$, $\mc{C}$ is precisely the $\mc{A}_0$ component in $\mc{A}$ in the splitting defined by $\wt{E}$.

\begin{lemm}
$\sigma(\wt{E}) = - \wt{E}$ and in the splitting defined by $\nabla$, the algebra $F_+$ is generated by $X + \rP(X)$ for sections $X$ of $T$, while $F_-$ is the span of elements of the form
\be
X - \rP(X) \ \ \textrm{and} \ \ \{X + \rP(X), Y - \rP(Y)\}.
\ee
Furthermore, $\rP$ follows all the properties of Theorem \ref{main:theo}.
\end{lemm}
\begin{lproof}
The proof will proceed by proving a series of interim results.

\begin{itemize}
\item $\sigma \wt{E} = - \wt{E}$.
\end{itemize}
Note that $\sigma$ sends the centre of $\mc{C} = \mc{A}_0$ to itself. Let $\pi$ be the projection $\mc{A} \to \mc{A}/\mc{A}_{(0)} \cong T$. By definition, $\pi \circ \sigma$ is bijective $T^* \to T$. Let $v_j$ be a section of $T^*_j$ in this given splitting. Then:
\be
j\sigma(v_j) = \sigma \{\wt{E},v_j\} = \{\sigma(\wt{E}), \sigma(v_j)\}.
\ee
Since $\sigma(\wt{E})$ must be a section of the centre of $\mc{A}_0$, it preserves the grading, implying that
\be
\pi(\{\sigma(\wt{E}), \sigma(v_j)\}) = \{\sigma(\wt{E}), \pi(\sigma(v_j))\}.
\ee
Thus the eigenvalues of $ad_{\sigma(\wt{E})}$ on $T$ are all strictly positive (as the eigenvalues of $ad_{\wt{E}}$ on $T^*$ are all strictly positive. Now let $w$ be a local section of $T^*_j$ that is an eigensection of $ad_{\sigma(\wt{E})}$. Then there exists a local eigensection $X$ of $T_{-j}$ such that $\{w,X\}$ is a nowhere zero section of $\mc{A}_0$. Since $\sigma(\wt{E})$ commutes with all of $\mc{A}_0$, the Jacobi identity gives:
\be
\{\{\sigma(\wt{E}),w\},X\}+\{w,\{\sigma(\wt{E}),X\}\} = 0.
\ee
Consequently the eigenvalues of $ad_{\sigma(\wt{E})}$ on $T^*$ are all strictly negative. Since $ad_{\sigma(\wt{E})}$ acts by multiplication by zero on $\A$ and $\sigma(v_j)$ is an eigensection of $ad_{\sigma(\wt{E})}$ with strictly positive eigenvalue, then $\sigma(v_j)$ must be a section of $T$ -- in other words, in this splitting, $\sigma(T^*) = \pi \circ \sigma(T^*)$ and thus $\sigma$ is a Lie algebra homomorphism $T^* \to T$.

By definition, $T_1^*$ generates the whole of $T^*$ by Lie algebra action. Consequently $\sigma(T_1^*)$ must be a generating bundle for $T$. Since the bracket preserves homogeneity and all sections of $T$ have strictly negative homogeneity, this means that the map
\be
T_1^* \to T_{-1}: v_1 \to (\sigma(v_1))_{-1},
\ee
must be a bijection. Thus
\be
\{\sigma(\wt{E}), (\sigma(v_1))_{-1} \} = \{\sigma(\wt{E}), (\sigma(v_1)) \}_{-1} = (\sigma \{\wt{E},v_1\})_{-1} = (\sigma(v_1))_{-1}.
\ee
Implying that $ad_{\sigma(\wt{E})}$ acts by multiplication by one on all of $T_{-1}$. This is the same action as that of $ad_{-\wt{E}}$; since $T_{-1}$ is a generating bundle for $T$, this means that the action of $ad_{\sigma(\wt{E})}$ and $ad_{-\wt{E}}$ match up on all of $T$, hence on all of $T^*$ and (trivially) on all of $\A_0$. So $ad_{\sigma(\wt{E})} = ad_{-\wt{E}}$, and since $\A$ is a semi-simple algebra bundle,
\be
\sigma(\wt{E}) = - \wt{E}.
\ee
This means that $\wt{E}$ is a section of $F_-$, and consequently that $ad_{\wt{E}}$ maps $F_-$ to $F_+$ and vice-versa. Thus $(ad_{\wt{E}})^2$ must maps $F_-$ and $F_+$ to themselves.

\begin{itemize}
\item $\rP$ is non-degenerate.
\end{itemize}
By the above, $\wt{E}$ is a section of $F_-$. This bundle must be preserved by $\onab$. Now if $\rP(X) = 0$ for a section $X$ of $T$, then in this splitting, $\onab_X \wt{E} = \{X,\wt{E}\} + \nabla_X \wt{E} = \{X,\wt{E}\} + \nabla_X \wt{E}$ is a section of $F_-$. Note that $\nabla_X \wt{E} = 0$, as $\nabla$ comes from a connection on a $G_0$ principal bundle. This implyes that $\{X,\wt{E}\}$ is a section of $F_-$. Consequently $\sigma$ must map a subbundle of $T$ to itself: impossible as $\sigma$ is a bijection between $T$ and $T^*$.

\begin{itemize}
\item Splitting $T$ into homogeneous components, if $X_{-a} \in \Gamma(T_{-a})$, then $\rP(X_{-i}, X_{-j}) = 0$ whenever $j \neq i$.
\end{itemize}
If the above statement fails, then there exist a section $X_{-j}$ of $T_{-j}$ and an $i \neq j$ such that $\rP(X_{-j})_i \neq 0$. If this is the case,
\be
\frac{1}{j^2} \{\wt{E}, \{\wt{E}, \onab_{X_{-j}} \wt{E} \} \}- \onab_{X_{-j}} \wt{E}
\ee
is a section of $F_-$, since $\onab_{X_{-j}} \wt{E}$ is a section of $F_-$ and $(ad_{\wt{E}})^2$ must map $F_-$ to itself. It is non-vanishing since it must have $i(\frac{i^2}{j^2} - 1) \rP(X_j)_i$ as the homogeneity $i$ component. However, its homogeneity $-j$ component is $j(1-\frac{j^2}{j^2})X_{-j} = 0$, and it is easy to see that all the other non-positive homogeneities vanish. This makes it into a non-vanishing section of $T^* \cap F_-$, again contradicting the fact that $\sigma$ is a bijection $T^* \to T$.

\begin{itemize}
\item $\sigma$ acts as $\rP: T \to T^*$, $\{ \rP(X), \rP(Y) \} = \rP(\{X,Y\})$ and $\rP$ is symmetric.
\end{itemize}
Note that the previous results imply that $\tau = X_{-j} -\rP(X_{-j}) = 1/j(\onab_{X_{-j}} \wt{E}$ is a section of $F_-$ and that $-\frac{1}{j} ad_{\wt{E}} \tau = X_{-j} + \rP(X_{-j})$ is a section of $F_+$. In other words, $F_+$ contains all sections of the form $X + \rP(X)$, while $F_-$ contains all sections of the form $X - \rP(X)$. This means that the map $\sigma : T \to T^*$ is given by $X \to \rP(X)$, and the fact that $\sigma$ is an algebra involution gives the relation $\{ \rP(X), \rP(Y) \} = \rP(\{X,Y\})$.

Now consider the Killing form $B$. By definition, $B(\tau, \nu) = $ trace $ad_{\tau} ad_{\nu}$. Since $\sigma$ preserves the Lie bracket, it must also preserve the invariant $B$, so $B(\sigma \tau, \sigma \nu) = B(\tau, \nu)$. Inserting $X$ and $Y$ into this and using $\sigma = \rP$ on $T$,
\be
\rP(X) \llcorner Y = B(\rP(X), Y) = B(X, \rP(Y)) = \rP(Y) \llcorner X.
\ee
So $\rP$ is symmetric.

\begin{itemize}
\item The algebra $F_+$ is generated by $X + \rP(X)$ for sections $X$ of $T$, while $F_-$ is the span of elements of the form $\{X + \rP(X), Y - \rP(Y)\}$.
\end{itemize}
$T \oplus T^*$ generates all of $\mc{A}$ by the Lie bracket, hence the result follows since we have fully defined the action of $\sigma$ on $T \oplus T^*$.

\begin{itemize}
\item $\rP$ defined an involution on $\mc{A}_0$ by conjugation on $T \otimes T^*$. This involution is the same as the restriction of $-\sigma$.
\end{itemize}
We know that $\sigma$ maps $\mc{C} = \mc{A}_0$ to itself. Then let $A$ be a section of $\mc{A}_0$; since $\sigma$ is an algebra involution,
\be
\{\sigma(A),X\} = \sigma(\{A,\sigma(X)\} = \rP(A(\rP(X))).
\ee

Note that this construction also works for the conjugation action of $\rP_j$ for any subbundles $T_{-j} \otimes T^*_j$ on which $\mc{A}_0$ acts faithfully -- $T_{-1} \otimes T^*_1$, for instance.

\begin{itemize}
\item $\nabla \rP = 0$, and the holonomy algebra bundle of $\nabla$ is contained in $\mc{B} = F_+ \cap \mc{A}_0$.
\end{itemize}
$\onab$ preserves $F_+$, and $X + \rP(X)$ evidently do so as well. Thus $\nabla_X$ must preserve $F_+$, implying that for a section $Y$ of $T$, $\rP(\nabla_X Y) = \nabla_X \rP(Y)$.

This implies that the holonomy algebra bundle of $\nabla$ is contained in $\mf{so}(\rP) \cap \mc{A}_0 = F_+ \cap \mc{A}_0 = \mc{B}$. This also means that $\nabla$ must preserve a volume form (in this instance, det $\rP$).
\end{lproof}

\subsection{Transitivity of Tractor bundles}
We take a pause now from Einstein involutions, to analyse some of the properties we will be needing later. For we will be generalising from the properties of $G$, while using as little as possible the properties of $P$. To do so, we need some universal properties of $\mc{T}$, not dependent on the choice of parabolic subalgebras. The most used will be the concept of cotransitivity. One immediate consequence of cotransitivity will be a restriction on the size of of the standard Tractor bundle $\mc{T}$.

\begin{defi}[Cotransitivity]
Let $\mf{g}$ be a Lie algebra, and $V$ a representation of dimension $m$. For an element $v\in V$, denote by $v^{\mf{g}}$ the orbit of $v$ under $\mf{g}$ -- it is a vector space, since $\mf{g}$ is. We say that $\mf{g}$ is $d$-cotransitive on $V$ if
\be
d = \textrm{ max}_{v\in V, v \neq 0} \{ \textrm{codimension of } v^{\mf{g}} \textrm{ in $V$} \}.
\ee
\end{defi}

\begin{prop}
Assume that $\mf{g}$ is $d$-cotransitive on $V$, and denote $r(V)$ the real dimension of $V$. If $V_l$ is the subbundle with highest homogeneity in $V$ (equivalently, the smallest subbundle in the natural filtration of $V$), then
\be
r(V) \leq r(\mf{g}_-) + r(V_l) + d.
\ee
\end{prop}
\begin{proof}
Given an element $v$ of $V_l$, the span $v^{\mf{g}}$ of $v$ under the action of $\mf{g}$ is of co-rank at most $d$ in $V$. However, $\mf{g}_{(0)}$ maps $V_l$ to itself; only the action of $\mf{g} / \mf{g}_{(0)} = \mf{g}_-$ can map $v$ non-trivially to $V / V_l$. The two inequalities
\be
r(V) &\leq& r(v^{\mf{g}}) + d \\
r(v^{\mf{g}}) &\leq& r(V_l) + r(\mf{g}_-),
\ee
then give the result.
\end{proof}

\begin{cor}
If we shift to the vector bundle point of view, with $\mc{V} = \mc{P} \otimes_P V$, $\mc{V}_l = \mc{P} \otimes_P V_l$ and (of course) $T = \mc{P} \otimes_P \mf{g}_-$, and where $r$ now denotes the real rank of a bundle, the preceding proposition implies that:
\be
r(\mc{V}) \leq r(T) + r(\mc{V}_l) + d.
\ee
\end{cor}

\begin{lemm} \label{co:trans}
The algebras $\mf{sl}(m,\mbb{F})$ and $\mf{sp}(2m, \mbb{F})$, $m>1$, are $0$-cotransitive (i.e. transitive) on their standard representations. The complex form of $\mf{so}$ is $2$-cotransitive on its standard representation, and $\mf{so}^*$ is $3$-cotransitive. All the other simple, non-exotic Lie algebras are $1$-cotransitive on their standard representations.
\end{lemm}
\begin{lproof}
This comes directly from differentiating the properties of the corresponding Lie groups. Both $SL(m,\mbb{F})$, $m>1$ and $Sp(2m,\mbb{F})$ are transitive on their standard representations, if we exclude the origin; differentiating this around any non-zero element yields the transitivity of their algebras.

All the other algebras ($\mf{so}(p,q)$, $\mf{so}(m,\mbb{C})$, $\mf{so}^*(2m)$, $\mf{su}(p,q)$, and $\mf{sp}(p,q)$) preserve a metric $h$ of some form on their standard representation $W$. And their corresponding groups all act transitively on connected components of the sets $W_{\lambda} = \{ w \in W | h(w,w) = \lambda \}$, again after excluding the origin (the proof of this in the $\mf{so}^*(2m)$ case will be detailed later; in the other cases the result is known, and with strictly simpler proofs).

There is an underlying real metric $h_{\mbb{R}} = Re(h)$. Except for $\mf{so}(m,\mbb{C})$ and $\mf{so}^*(2m)$, $h(w,w)$ is always real, so $h_{\mbb{R}}(w,w) = h(w,w)$. Consequently, differentiating this local transitivity on $W_{\lambda}$, we get the result that $\mf{g}$ will map $w \neq 0$ to all of $w^{\perp}$, where $\perp$ is taken with respect to $h_{\mbb{R}}$. Hence these algebras are $1$-cotransitive on their standard representations.

For $\mf{so}(m,\mbb{C})$ the same result applies, except that we need to use $h$ to define $\perp$ rather than $h_{\mbb{R}}$. Consequently $\mf{so}(m,\mbb{C})$ is $2$-cotransitive.

The algebra $\mf{so}^*(2m)$ is a bit more subtle. It can be seen as the algebra that acts on a quaternionic space $V \cong \mbb{H}^m$, preserving the quaternionic multiplication $i$, $j$ and $k$, as well as a real metric $h$ that is hermitian with respect to $j$ and symmetric with respect to $i$ and $k$ (neither $h$, nor the choices of $i$, $j$ and $k$ are canonical). Fix a orthonormal basis of the form $\{e_1, ie_1, je_1, ke_1, e_2, ie_2, \ldots, ke_n\}$, and use the multiplication by $i$ to set an isomorphism $V \cong \mbb{C}^{2m}$. Then this algebra is expressed, in matrix form, as:
\beqa \label{so:matrix}
\left( \begin{array}{cccc}
A_{11} & A_{12} & \ldots & A_{1m} \\
A_{21} & A_{22} & \ldots & A_{2m} \\
\vdots & \vdots & \ddots & \vdots \\
A_{m1} & A_{m2} & \ldots & A_{mm}
\end{array}
\right),
\eeqa
with
\be
A_{ll} = \left( \begin{array}{cc}
0 & -a_{ll} \\
a_{ll} & 0 \\
\end{array}
\right) \ \ \textrm{and} \ \ A_{lm} = - A_{ml}^t = \left( \begin{array}{cc}
\alpha_{lm} & -\overline{\beta_{lm}}\\
\beta_{lm} & \overline{\alpha_{lm}}
\end{array}
 \right) \in \mf{csu}(2) \cong \mbb{H},
\ee
for real numbers $a_{ll}$ and complex numbers $\alpha_{lm}$ and $\beta_{lm}$. To get a more invariant definition, let us define $\tilde{h}$, a non-degenerate bilinear form on $V$ with values in $\mbb{H}$. It is defined as:
\be
\tilde{h}(v,w) = -h(v,jw) + jh(v,w) +ih(v,kw) -kh(v,iw).
\ee
This is hermitian, as $\tilde{h}(av,bw) = a\tilde{h}(v,w)\overline{b}$ for any quaternions $a$ and $b$. Moreover, $\tilde{h}(v,w) = - \overline{\tilde{h}(w,v)}$ (notice the contrast with the standard quaternionic-hermitian metric, where that relationship would be $g(v,w) = \overline{g(w,v)}$). This ensures that $Im(\tilde{h})$ is symmetric, while $Re(\tilde{h})$ is skew -- thus $\tilde{h}(v,v)$, the $\tilde{h}$ norm-squared of $v$, is imaginary for all $v$ in $V$. The above properties ensure that for any $v\in V$, $v\neq 0$, then the map
\be
V \to Im(\mbb{H}) , w \to \tilde{h}(w,v),
\ee
is surjective.

Now let $u$ be an element of $V$, and assume that $\tilde{h}(u,u) = z \neq 0$. By real scaling of $u$, we may assume that $z$ is of unit norm (hence $z^{-1} = -z$). We may now define a new real metric by replacing $h$ with
\be
h'(v,w) = \frac{1}{2z}(\tilde{h}(v,w) + \tilde{h}(zv,zw)) = -Re (z \cdot \tilde{h}(v,w)).
\ee
Replace $i$, $j$ and $k$ with $j' = z$ and $i'$ and $k'$ any unit imaginary quaternions that anti-commute with $z$ and each other (hence that are orthogonal to $z$ and each other in the standard norm). Under these conditions, $u, i'u, j'u$ and $k'u$ are orthonormal, and can be extended to an orthonormal basis of $V$ such that $\mf{so}^*(2m)$ is of the form detailed in equation (\ref{so:matrix}). In this basis, $u = (1,0,\ldots)^t$. This demonstrates that the span of $u$ under $\mf{so}^*(2m)$ maps $u$ onto $u^{\perp}$, with $\perp$ being taken via $\tilde{h}$.

Now assume $\tilde{h}(u,u) = 0$, $u \neq 0$. Let $u'$ be any other null vector in $V$ that is not orthogonal to $u$. Then the vectors $v = u + u'$ and $w=u-u'$ are orthogonal, and have the property that $\tilde{h}(v,v) = -\tilde{h}(w,w) = 2\tilde{h}(u,u')$. By scaling $u$, we may ensure that $v$ and $w$ are of unit norm. Then we may set $z = \tilde{h}(v,v)$ and use the same procedure as previously to get the orthonormal set $\{v, i'v, j'v, k'v, (-i'w),$ $i'(-i'w), j'(-i'w), k'(-i'w)\}$. Extending this to a suitable orthonormal basis of $V$, as before, we now have $u = (1,0,0,0,0,1,0,\ldots)^t$.

Inserting this into the matrix form of $\mf{so}^*$, and a little work, reveals that this algebra maps $u$ onto $u^{\perp}$.

Since $u^{\perp}$ is of co-dimension three, we are done. Note that by integrating the above action, we can see that the group $SO^*(2n)$ is transitive on connected components with constant $\tilde{h}$-norm.
\end{lproof}

The remaining non-exotic simple algebras with are not $0$-cotransitive will be called \emph{metric} (for obvious reasons).

By these preceding results, we can affirm:
\begin{prop}
For $\mf{g} = \mf{sl}(m,\mbb{F}), \mf{sp}(2m,\mbb{F})$, $\mbb{F} = \mbb{R}, \mbb{C}$ and $\mf{g} = \mf{sl}(m, \mbb{H})$, $m>1$,
\be
r(\mc{T}) \leq r(T) + r(H_l).
\ee
For all other non-exotic simple $\mf{g}$, apart from $\mf{so}(h,\mbb{C})$ and $\mf{so}^*(2m)$,
\be
r(\mc{T}) \leq r(T) + r(H_l) + 1
\ee
for $\mf{so}(h,\mbb{C})$
\be
r(\mc{T}) \leq r(T) + r(H_l) + 2,
\ee
and for $\mf{so}^*(2m)$
\be
r(\mc{T}) \leq r(T) + r(H_l) + 3.
\ee
\end{prop}

\subsection{Existence of Einstein involutions for certain holonomy reductions}
There are several possible involutions that can be defined on the bundle $\mc{A}$; for instance, conjugation by a complex structure. Some preserved structures, however, generate Einstein involutions in a natural way. To define them (especially for metric algebras), we shall need the technical concept of the \emph{image degree}:
\begin{defi}
Let $H_{(k)}$ be the filtration component of minimal homogeneity $k$ of a Tractor bundle $\mc{V}$ (see the end of section \ref{parabolic:section}). Then $H_{(k)}$ is called an \emph{image bundle} of $T^* \subset \mc{A}$ if for all local never-zero sections $v$ of $T^*$, there exists local $\xi\in\Gamma(\mc{V})$ such that $v \cdot \xi$ is a never-zero local section of $H_{(k)}$. The full bundle $\mc{V}$ is always an image bundle as $\mc{A}$ is a simple algebra bundle, so has a faithful action on all bundles associated to it. This shows that image bundles exist for every Tractor bundle.

The \emph{image degree} of $\mc{V}$ is the (unique) $r$ such that $H_{(r)}$ is an image bundle of $T^*$ while $H_{(r+1)}$ is not (recall that these filtered subbundles are indexed by the lowest eigenvalue of a given grading section $\wt{E}$ acting on them; this number does not depend on the choice of grading section).
\end{defi}
\begin{exa} In the projective case, the subbundle $\mc{E}[\mu] \subset \mc{T}$ (see \cite{mepro1}; $\mu = \frac{n}{n+1}$) is an image subbundle, as the action $T^* \cdot \T \to \mc{E}[\mu]$ is given by the trace of $T^*$ with $T$, evidently surjective for all sections of $T^*$. Since $\mc{E}[\mu]$ is the filtered subbundle of highest homogeneity, $-\mu$ is the image degree for the standard Tractor bundle in projective geometry (the fact that it is $-\mu$ not $\mu$ come from conflicting conventions: the actions of a grading section on $\mc{E}[\lambda]$ is via multiplication by $-\lambda$).
\end{exa}
\begin{exa}
In conformal geometry, let $\mc{E}[-1] \subset \T$ be the highest homogeneity bundle. The map $T^* \cdot (T[-1] \rtimes \mc{E}[-1]) \to \mc{E}[-1]$ is given by contracting $T^*$ with $T$, and is thus surjective for all sections of $T^*$. Again, because of the sign convention, the image degree of $\T$ for conformal geometry is $+1$.
\end{exa}

\begin{lemm} \label{ima:met}
If the algebra bundle $\mc{A}$ preserves a metric $h$ on $\mc{V}$, then the image degree of $\mc{V}$ is strictly positive.
\end{lemm}
\begin{lproof}
Pick a Weyl structure $\nabla$, and a consequent splitting of $\mc{V}$ and $\mc{A}$. Let $v$ be a never-zero local section of $T^*_k \subset \mc{A}$. Since $\mc{A}$ is a simple algebra bundle, $v$ has a non-trivial action on $\mc{V}$. Consequently there exists a homogeneous local section $\xi$ of $H_j \subset \mc{V}$ such that $v \cdot \xi$ is never-zero, and of homogeneity $j+k$. Since $h$ has homogeneity zero and is non-degenerate, there exists a local section $\eta$ of homogeneity $-(j+k)$ such that
\be
1 = h(v \cdot \xi, \eta) = h(\xi, v \cdot \eta).
\ee 
Consequently $v \cdot \eta$ is a never-zero local section of $H_{-j}$. Since $k>0$, one at least of $j+k$ and $-j$ is strictly positive, so one of $v \cdot \xi$ and $v \cdot \eta$ is of strictly positive homogeneity. As the homogeneity degrees of $\mc{V}$ are discreet (indeed, their difference are integers), this implies that the image degree of $\mc{V}$ is strictly positive.
\end{lproof}

We are now able to phrase the major result on creating Einstein involutions in the simple, non-exceptional cases:
\begin{theo} \label{theo:existance}
Let $\mc{T}$ be the standard Tractor bundle coming from the standard representation of $G$.

For $\mf{g} = \mf{sl}(m, \mbb{R})$, and $\mf{sl}(m,\mbb{C})$ any non-degenerate metric $g$ that lies in an irreducible component of $\odot^2 \mc{T}$, and that has $\onab g = 0$, generates an Einstein involution on an open, dense subset of $M$.

For $\mf{g} = \mf{sl}(m, \mbb{H})$, assume $g$ non-degenerate. The bundle $\mc{T} \odot \mc{T}$ splits into two irreducible components, one isomorphic at each point to $(\mbb{H}^m \otimes_{\mbb{H}} \overline{\mbb{H}}^m) \cap ({\mbb{H}}^m \odot {\mbb{H}}^m)$. If $g$ lies in this component, it is hermitian with respect to all the complex structures of the quaternionic structure, and the same result hold as above. If $g$ lies in the other irreducible component, then generically we can construct another metric $g$ such that there exists complex structures $I,J$ and $K$, compatible with the quaternionic structure and obeying the quaternionic identities, such that $g$ is hermitian with respect to $I$ and symmetric with respect to $J$ and $K$. Then the same result holds.

For $\mf{g} = \mf{sp}(2m, \mbb{F})$, let $J$ be a complex structure on $\mc{T}$ such that the natural skew-form $w$ is hermitian with respect to it, and such that $\onab J = 0$. If $\mbb{F} = \mbb{C}$, then we further require that $J$ be complex-linear or complex-hermitian with respect to multiplication by $i$. Then $J$ generates an Einstein involution on an open, dense subset of $M$.

For $\mf{g}$ any of the other non-exceptional simple Lie algebras (which are $\mf{so}(p,q)$, $\mf{so}(m,\mbb{C})$, $\mf{so}^*(2m)$, $\mf{su}(p,q)$ and $\mf{sp}(p,q)$), let $K$ be non-degenerate subbundle $K \subset \mc{T}$, that is preserved $\onab$ and preserved by any complex structures associated to $\mf{g}$. If the rank and co-rank of $K$ in $\mc{T}$ are both higher that the rank of $\T_{(r)}$, for $r$ the image rank of $\mc{T}$, then $K$ defines an Einstein involution on an open, dense subset of $M$.
\end{theo}
This section will be devoted to proving that.

\begin{itemize}
\item $\mf{g} = \mf{sl}(m, \mbb{F})$
\end{itemize}
If $\mbb{F} = \mbb{R}$, then conjugation by $g$ defines an involution of $\mc{A} = \mc{T} \otimes_0 \mc{T}^*$. If $\mbb{F} = \mbb{C}$, then $g$ is either $J$-symmetric or $J$-hermitian, since it is in an irreducible component of $\odot^2 \mc{T}$. If $\Psi$ is a section of $\mc{A}$, then
\be
(g \Psi g^{-1}) J = \pm 1 (g \Psi Jg^{-1}) = \pm 1 (g J \Psi g^{-1}) = (\pm 1)^2 J(g \Psi g^{-1}).
\ee
This implies that $g \Psi g^{-1}$ is also a section of $\mc{A}$, so conjugation by $g$ does preserve $\mc{A}$.

If $\mbb{F} = \mbb{H}$, then $\odot^2 \T^*$ splits into two irreducible bundles. The smallest one is $\mc{T} \odot_{\mbb{H}} \overline{\mc{T}}$, which is the bundle of totally hermitian metrics. By the above argument, the requirements we have put on $I,J$ and $K$ guarantee that conjugation by $g$ preserves all these complex structures, hence preserves $\A$.

\begin{rem}
Note that the requirements that be there exist complex structures $I,J$ and $K$, compatible with the quaternionic structure and obeying the quaternionic identities, such that $g$ is hermitian with respect to $I$ and symmetric with respect to $J$ and $K$, is not a strong one. Given any $g$, $I$, $J$ and $K$, we can form the projection
\be
g(v,w) \to g(v,w) + g(Iv,Iw) -g(Jv,Jw) - g(Kv,Kw),
\ee
which has the desired properties. The projection in non-degenerate on an open dense set of this irreducible component, and anything that preserves $g$ and the complex structures will preserve this projection.
\end{rem}

Let $\sigma$ be this involution generated by conjugation, and we now appeal to the following lemma:
\begin{lemm} \label{ein:evo}
If a metric $g$ is non-degenerate on $\T_{(k)}$ for all $k \geq r$ for $r$ the image degree of $\mc{T}$, then $\sigma$ is an Einstein involution.
\end{lemm}.
\begin{lproof}
We merely need to show that $\sigma(T^*) \cap \mc{A}_{(0)} = 0$. This will be proved by contradiction. Let $\nu$ be non-zero section of $T^*$ such that $A = \sigma(\nu)$ is a section of $\mc{A}_{(0)}$. We will demonstrate that the image of $\mc{T}$ under $\nu$ does not intersect $\T_{(r)}$, contradicting the definition of the image degree. The proof will proceed by induction down the homogeneity degrees of the bundles $\T_{(j)} \subset \mc{T}$. Let $l$ be the maximal homogeneity degree of $\mc{T}$.

Let $P_j$ be the proposition that:
\begin{enumerate}
\item the image of $\mc{T}$ under $\nu$ does not intersect $\T_{(l-j)}$, and
\item $A$ has a trivial action on $\T_{(l-j)}$.
\end{enumerate}
Note that since $\nu$ has strictly positive homogeneity, the first statement implies that $\nu$ has trivial action on $\T_{(l-j)}$. We need to show that if $l-j>r$, $P_j$ implies $P_{j+1}$.

Let $u$ and $v$ be sections of $\T_{(l-j-1)}$ and $t$ be any section $\T_{(l-j)}$. The relation
\be
g(\nu \cdot u, t) = g(u, A \cdot t) = 0,
\ee
demonstrates that $\nu$ has a trivial action on $\T_{(l-j-1)}$ (since $\nu$ must increase homogeneity, it must map sections of $\T_{(l-j-1)}$ into $\T_{(l-j)}$). Then the relation
\be
g(u, A \cdot v) = g(\nu \cdot u, v) = 0,
\ee
shows that $A$ itself must have trivial action on $\T_{(l-j-1)}$. Now let $\eta$ be any section of $\mc{T}$; the relation
\be
g(\nu \cdot \eta, v) = g(\eta, A \cdot v) = 0,
\ee
gives the final condition that the image of $\mc{T}$ under $\nu$ does not intersect $\T_{(l-j-1)}$.

The above proof also demonstrate that $P_{-1}$ implies $P_0$, thus completing the induction, since $P_{-1}$ is trivially true as $\T_{(l+1)} = 0$.
\end{lproof}

Now, to demonstrate that we have an Einstein involution on an open dense subset of $M$, we require:
\begin{lemm}
Let $g$ be any metric on $\T$ that is non-degenerate, preserved by $\onab$ for a parabolic geometry with $\mf{g} = \mf{sl}(m,\mbb{F})$. Then on an open dense subset of $M$, $g$ is also non-degenerate on all subbundles $\T_{(k)}$ of the filtration of $\T$.
\end{lemm}
\begin{lproof}
Let $M_k$ be the subset of $M$ where $g$ is non-degenerate on $\T_{(k)}$. Since non-degeneracy is an open condition, $M_k$ is open.

Reasoning by contradiction, let $x \in M - M_k$, and let $\eta$ be a local non-vanishing section of $\T_{(k)}$ such that $g(\eta,\tau) = 0$ for all sections $\tau$ of $\T_{(k)}$ in a neighbourhood of $x$. Then for any section $X$ of $T$:
\be
0 = \nabla_X g(\eta, \eta) = 2g(\nabla_X \eta, \eta).
\ee
Since $\eta$ is orthogonal to $\T_{(k)}$, $g(-,\eta)$ descends to a section of $\mc{T}/\T_{(k)}$. In that setting, $g(\nabla_X \eta, \eta)$ descends to $g(X \cdot \eta, \eta)$, for the well defined action $T = \mc{A}/\mc{A}_{(0)}: \T_{(k)} \to \mc{T} / \T_{(k)}$. However, since $\mc{A}$ is transitive on $\mc{T}$ (as $\mf{sl}(m,\mbb{F})$ is $0$-cotransitive, and since $\mc{A}_{(0)}$ maps $\T_{(k)}$ to itself, $T \otimes \eta$ must map surjectively to $\mc{T} / \T_{(k)}$. Consequently, $g(-,\eta) = 0$, which is contradicted by the non-degeneracy of $g$. This implies that $M_k$ must be an open dense subset of $M$.
\end{lproof}

Thus the conditions of Lemma \ref{ein:evo} hold on the open dense set $\cap_{k = r}^l M_k$.

\begin{itemize}
\item $\mf{g} =\mf{sp}(2n,\mbb{C})$.
\end{itemize}
Let $w$ be the underlying real skew-form on $\mc{T}$. Since $J$ anti-commutes with $w$, $g = w \circ J$ is a symmetric bilinear form on $\mc{T}$. Since $w$ is non-degenerate and $J$ is an automorphism, $g$ is a metric. Furthermore,
\be
w^{-1} g = -g^{-1} w,
\ee
so conjugation by $w$ and $g$ commute. Since furthermore $g$ must be complex linear or complex hermitian, conjugation by $g$ must be complex linear, implying that conjugation by $g$ is an involution of $\mc{A}$.

Then since the action of $\mc{A}$ on $\mc{T}$ is still transitive in this case, the proof proceeds exactly as in the special linear case.

\begin{rem}
For $\mbb{F} =\mbb{C}$, if $J$ is complex hermitian, the holonomy algebra of $\onab$ reduces to $\mf{sp}(p,q)$. If it is complex linear, the holonomy algebra of $\onab$ reduces to $\mf{gl}(n,\mbb{C})$.
\end{rem}

\begin{itemize}
\item $\mf{g}$ is one of the metric non-exceptional algebras, preserving an underlying real metric $h$ on $\mc{T}$.
\end{itemize}
These $\mf{g}$ are $\mf{so}(p,q)$, $\mf{so}(m,\mbb{C})$, $\mf{so}^*(2m)$, $\mf{su}(p,q)$, and $\mf{sp}(p,q)$. In all these cases, $\mc{A}$ is defined as the subbundle of $\mc{T}\otimes\mc{T}^*$ commuting with conjugation by $h$ and with the complex structures.

We will use the bundle $K$ to generate another metric $g$ such that conjugation by $g$ commutes with conjugation by $h$ and with the complex structures; this is enough to demonstrate that conjugation by $g$ preserves $\mc{A}$ and thus define an involution $\sigma$ on it.

Let $K^{\perp}$ be the orthogonal complement of $K$ using the bilinear form $h$; the non-degeneracy of $K$ implies that $K \oplus K^{\perp} = \mc{T}$ and that $K^{\perp}$ is preserved by the complex structures. The bilinear form splits as $h_1 + h_2$, $h_1$ a section of $\odot^2 K$, and $h_2$ a section of $\odot^2 K^{\perp}$. Since $\onab$ must preserve $K$ and $K^{\perp}$ and $\onab h = 0$,
\be
\onab h_1 = \onab h_2 = 0.
\ee
And since the $K$ are preserved by the complex structures, each $h_j$ must be hermitian or symmetric with respect to the complex structures.

Then define $g = h_1 - h_2$. Since $g$ is hermitian or symmetric with the complex structures, conjugation by $g$ on $\mc{T}\otimes\mc{T}^*$ commutes with complex multiplication. Moreover:
\be
h g^{-1} = (h_1 + h_2)(h_1^{-1} - h_2^{-1}) = Id_K - Id_{K^{\perp}} = g h^{-1},
\ee
so conjugation by $h$ and $g$ commute.

Lemma \ref{ein:evo} still applies for this $g$. Since we have the natural bilinear form $h$, Lemma \ref{ima:met} implies that the image degree $r$ of $\mc{T}$ is positive. We shall prove that if $g$ cannot degenerate on $H_{(k)}$, for $k \geq r>0$.

\begin{lemm}
Assume the metric $g$ is degenerate on $H_{(j)}$ for any $j>0$. Let $\xi$ be a section of $H_{(j)}$ on which $g$ degenerates. Then $\xi$ is a section of $K$ or $K^{\perp}$.
\end{lemm}
\begin{lproof}
By definition,
\be
0 = \onab_X g(\xi,\xi) = 2 g (\onab_X \xi, \xi).
\ee
The map $g(-,\xi)$ descends, by the degeneracy of $\xi$, to a section of $(\mc{T} / H_{(j)})^*$. So only the $\mc{T} / H_{(j)}$ component of $\onab_X \xi$ matters in the previous equation. This is simply $X$ acting algebraically on $\xi$. Only the $T$ component of $\mc{A}$ reduces homogeneity, and $\mc{A}$ maps $\xi$ onto $\xi^{\perp}$ by Lemma \ref{co:trans}.

Here, we have used $\tilde{h} = h$ to define $\perp$, apart from the algebras $\mf{so}(m,\mbb{C})$, where we have used the complex metric $\tilde{h}(v,w) = h(v,w) - ih(v,iw)$, and $\mf{so}^*(2m)$, where we have used the quaternionic `metric' $\tilde{h}(v,w) = -h(v,jw) + jh(v,w) +ih(v,kw) -kh(v,iw)$.

This means that as $X$ varies, $\onab_X \xi/H_{(j)}$ must span $\xi^{\perp} / H_{(j)}$. Consequently, $g(\eta, \xi) = 0$ whenever $\tilde{h}(\eta, \xi) = 0$. If $\tilde{h} = h$, this implies that $g(-,\xi)$ is a multiple of $\tilde{h}(-,\xi)$. Hence $\xi$ is an eigensection of $h^{-1} g$, thus a section of $K$ or $K^{\perp}$.

For $\mf{so}(m,\mbb{C})$, we may define $\tilde{g}(v,w) = g(v,w) - ig(v,iw)$. Now the complex structure is preserved by $\mc{A}$, hence multiplication by $i$ must be of homogeneity zero, implying that $i$ is an isomorphism of $H_{(j)}$ (by definition it is an isomorphism of $K$ and $K^{\perp}$). In particular $g$ must degenerate on $i\xi$, and thus $g(\eta,i\xi) =0$ whenever $\tilde{h}(\eta,i\xi) = i\tilde{h}(\eta,\xi) = 0$. This implies that $\tilde{g}(-,\xi)$ is a complex multiple of $\tilde{h}(-,\xi) = 0$, and thus that $g(-,\xi)$ is a real multiple of $h(-,\xi) = Re(\tilde{h}(-,\xi))$. Thus, as above, $\xi$ is a section of $K$ or $K^{\perp}$.

For $\mf{so}^*(2n)$, we may similarly define $\tilde{g}(v,w) = -g(v,jw) + jg(v,w) +ig(v,kw) -kg(v,iw)$. Now $g$ is $i$ and $k$ linear, and $j$ hermitian, so $g$ must degenerate on $i\xi$, $j\xi$ and $k\xi$ on $H_{(j)}$. Thus $\tilde{g}$ degenerates on $\xi$. Differentiating $\tilde{g}(\xi,\xi)$ implies that $\tilde{g}(\eta,\xi) = 0$ whenever $\tilde{h}(\eta,\xi) = 0$. Since $\tilde{g} = \tilde{h}_1 - \tilde{h}_2$ (using the evident definitions of $\tilde{h}_1$ and $\tilde{h}_2$), what we have is the relation
\be
\big(\tilde{h}_1(\eta,\xi) + \tilde{h}_2(\eta,\xi) = 0\big) \Rightarrow \big( \tilde{h}_1(\eta,\xi) = 0 \ \ \textrm{and} \ \ \tilde{h}_2(\eta,\xi) = 0 \big).
\ee
However, if $\xi$ is not a section of $K$ or $K^{\perp}$ (thus $\xi_1 \neq 0$ and $\xi_2 \neq 0$), fix any $\eta_1$ in $K$ such that $\tilde{h}_1(\eta_1,\xi) = \tilde{h}_1(\eta_1,\xi_1)= i$. We may do this as $\tilde{h}_1$ is non-degenerate on $K$ and $\tilde{h}_1(-,\xi_1)$ maps surjectively onto $\mbb{H}$ (see the proof of Lemma \ref{co:trans}). Then the equation $\tilde{h}_2(\eta_2,\xi_2) = -i$ must have a solution for a certain $\eta_2$ of $K^{\perp}$. Then $\eta = \eta_1 + \eta_2$ gives us a violation of the above condition, hence a contradiction. Thus $\xi$ must be a section of $K$ or $K^{\perp}$.
\end{lproof}

From the previous Lemma, it suffices to show that $K$ and $K^{\perp}$ do not intersect $H_{(r)}$. The conditions on the rank and co-rank of $K$, demonstrate that this is equivalent with them being transverse on an open dense subset of $M$. The next Lemma establishes that fact.

\begin{lemm}
Let $K$ be any preserved subset of $\mc{T}$, for any parabolic such that $\mc{A}$ fixes a (real, complex or quaternionic) metric $\tilde{h}$ on $\mc{T}$. Then on an open dense subset of $M$, $K$ and $H_{(j)}$ are transverse, for \emph{all} $j$.
\end{lemm}
\begin{lproof}
Let $\xi$ be a never-zero local section of $H_{(j)}$. For any section $X$ of $T$, the operator $\onab_X$ operates on $\xi$. The only component of $\onab_X$ that maps $\xi$ non-trivially into sections of $\mc{T} / H_{(j)}$ is $X$ itself, acting algebraically. Since $T$ is the only component of $\mc{A}$ that reduces homogeneity, and since the action of $\mc{A}$ on $\mc{T}$ maps any section to its $\tilde{h}$-orthogonal complement, the span of $\xi$ and its Tractor derivatives, quotiented out by $H_{(j)}$, is
\be
\xi^{\perp} / H_{(j)}
\ee
Since $\tilde{h}$ is of homogeneity zero, if $j \leq 0$, this is just $\mc{T} / H_{(j)}$. If not, then the span of $\xi$ under sections $X$ of minimal homogeneity $\geq -j$ is $H_{(0)} / H_{(j)}$. The span of this under a second iteration of $\onab_X$ is then $\mc{T} / H_{(j)}$.

So, the span under iterated Tractor derivatives for any non-zero section $\xi$ of $H_{(j)}$, projects surjectively onto $\mc{T} / H_{(j)}$.

Now assume that $K$ and $H_{(j)}$ are not transverse on some open set of $M$. This means that there exists a local nowhere-zero section $\xi$ of $K \cap H_{(j)}$ and that the projection of $K$ to $\mc{T} / H_{(j)}$ is not surjective. But this is contradicted by the fact that iterated derivatives of $\xi$ must be sections of $K$, as $\onab$ preserves $K$.
\end{lproof}

\begin{rem}
The results of this section give the results for holonomy reduction cited in tables \ref{table:one} and \ref{table:two}.
\end{rem}

\subsubsection{Examples of Einstein involutions for various geometries}\label{invo:exam}

The normality condition is detailed in both \cite{CartEquiv} and \cite{TCPG}; it means that the curvature $\ora{R}$ of $\onab$ is closed under the algebraic Lie algebra co-differential $\partial^*: \wedge^2 T^* \otimes \mc{A} \to T^* \otimes \mc{A}$. There is a unique regular normal Tractor connection for each parabolic geometry, making it a uniqueness condition similar to torsion-freeness for a metric connection. In this section, we will seek to build examples of normal Tractor connections preserving an Einstein involution.

This is the case, of course, for conformal and projective geometries: \cite{mepro1} demonstrates that there exists (on an open, dense set) an Einstein connection in the projective class of a projective manifold if and only if the Tractor connection preserves a metric on $\mc{T}$. Here we define an Einstein connection as a connection $\nabla$ with $\mathsf{Ric}^{\nabla}$ non-degenerate and
\be
\nabla \mathsf{Ric}^{\nabla} = 0.
\ee
This makes $\nabla$ into an Einstein connection for the metric $g = \mathsf{Ric}^{\nabla}$. The metric on $\T$ the generates the Einstein involution by Theorem \ref{theo:existance}.

In the conformal class, there are of course standard Einstein metric generating a preserved Tractor, and hence an extra preserved metric on $\mc{T}$. But the Einstein conformal product decomposition \cite{mecon} also generates preserved bundles and hence an extra preserved metric on $\mc{T}$, and thus an Einstein involution. The conditions on the magnitude of the Einstein coefficients in paper \cite{mecon} can be naturally interpreted as equivalent with the condition that $\rP$ defines the involution $\sigma$ on $\mc{A}_0$ by minus conjugation.

Normal contact-projective geometries are equivalent with normal projective geometries such that $\onab$ preserves a skew-form $\nu$ on the $\mc{T}$ (see \cite{fox} and also \cite{mepro1}). Thus if the projective holonomy algebra of $\onab$ reduces to $\mf{su}(p,q)$ (which is the case for projectively Sasaki-Einstein manifolds), then it generates an Einstein involution for a contact-projective geometry. This is particularly significant, as this is an Einstein involution on a $|2|$-graded geometry.

For CR geometries ($\mf{g} = \mf{su}(p+1,q+1), \mf{p} = \mf{cu}(p,q) \rtimes \mbb{C}^{p,q} \rtimes \mbb{R}$), the standard Tractor bundle splits as $\mbb{C} \oplus \mbb{C}^{p,q} \oplus \mbb{C}^*$. The bundle of strictly positive homogeneity is of complex rank one, so must be the image bundle. Let $\tau$ be a non-isotropic section of $\mc{T}$ such that $\onab \tau = 0$. This implies that $\onab J \tau = 0$ and that the space $K$ spanned by $\tau$ and $J \tau$ is non-degenerate. Because of its rank and co-rank, Theorem \ref{theo:existance} gives us an Einstein involution on an open dense set of $M$. Papers \cite{complexconf}, \cite{govfef}, \cite{leitnersu} and \cite{leitnercom} demonstrate the properties and existence of geometries with such preserved Tractors. Again, this is an Einstein involution on a $|2|$-graded geometry.

An Einstein involution can be also defined for normal almost Grassmannian geometries of degree two, and their almost quaternionic analogues. These are the geometries given in terms of crossed nodes \cite{two} as $\grasstwogeo$. Now, paper \cite{CartEquiv} demonstrates that the lowest homogeneity component of the curvature of a normal Tractor connection must be harmonic: i.e. closed under $\partial^*$ and its dual $\partial$. By Kostant's version of the Bott-Borel-Weil theorem \cite{Kostant}, for almost Grassmannian geometries of degree two, the harmonic piece of $\wedge^2 T^* \otimes \mc{A}$ has two components, one in $\wedge^2 T^* \otimes T$, the second in $\wedge^2 T^* \otimes \mc{A}_{(0)}$. If we have no curvature in the first component, $\onab$ is a torsion-free Tractor connection \cite{TCPG}. In particular, since this geometry it is $|1|$-graded and non-projective, the preferred connections are precisely those torsion-free connection with structure bundle $\mc{G}_0$. We will construct the Einstein involution using the following theorem:
\begin{prop}
Let $(M, \mu)$ be a (pseudo-)Riemannian manifold with metric $\mu$ and holonomy algebra contained in $\mf{h} = \mf{sp}(p,q) \cdot \mf{sp}(1)$ -- a quaternionic Kh\"aler manifold. Then $M$ defines a torsion-free almost quaternionic geometry of degree two with an Einstein involution on the normal Tractor connection $\onab$.

Moreover, there exist manifold with such properties, such that $\onab$ is non-flat.
\end{prop}
\begin{proof}
Let $H$ be the Lie group $Sp(p,q).Sp(1)$ and let $\nabla$ be the Levi-Civita connection of $\mu$. Because of its holonomy algebra, $\nabla$ must preserve both bundles in the tensor product
\be
T = \mc{L} \otimes_{\mbb{H}} U,
\ee
where $\mc{L}$ is a left-quaternionic line bundle and $U$ a right-quaternionic bundle of same real dimension as $T$.

Let $\mc{H}$ be the frame bundle for $\nabla$, and define $\mc{A}_0 = \mc{H} \times_H \mf{g}_0$ and, via the Weyl structure $\nabla$,
\be
\mc{A} = T \oplus \mc{A}_0 \oplus T^* = \mc{H} \times_H \mf{g}.
\ee
The differential $\partial$ is well-defined as a map $T^* \otimes \mc{A} \to  \wedge^2 T^* \otimes \mc{A}$; homogeneity considerations means that it restricts to a map $T^* \otimes T^* \to \wedge^2 T^* \otimes \mc{A}$.

The composition $\partial^* \circ \partial$ is bijective as a map from $T^* \otimes T^*$ to itself. If $R$ is the curvature tensor of $\nabla$, then $\partial^* R \in \Gamma(T^* \otimes T^*)$, so we define $\rP$ to be the section of $T^* \otimes T^*$ such that $\rP = -(\partial^* \circ \partial)^{-1}(\partial R)$. This implies that
\be
\partial^* (R + \partial \rP) = 0.
\ee
Now $\partial R$ is just the Ricci-trace of $R$ (see \cite{TCPG} in the $|1|$-graded case). All quaternionic Kh\"aler manifolds are Einstein (see \cite{spein} \cite{berger}), so $\partial R$ is a multiple of $\mu$. Since $\partial^* \circ \partial$ are $G_0$-module isomorphisms, this implies that $\rP$ must also be a multiple of $\mu$, as this is the only irreducible line subbundle of $T^* \otimes T^*$.

Consequently $\nabla \rP = 0$. We then define a Tractor connection on $\mc{A}$ as
\be
\onab_X = X + \nabla_X + \rP(X).
\ee
The curvature of this connection is $(0, R + \partial \rP, 0)$. This is $\partial^*$-closed, so $\onab$ is normal. By the uniqueness result for normal Cartan connections \cite{CartEquiv}, this is the unique normal Cartan connection for this geometry.

We now aim to show that there is an Einstein involution on $\mc{A}$, generated by $\rP: T \to T^*$, $\rP^{-1}: T^* \to T$ and by the conjugation action of $\rP$ on $\mc{A} \subset T^* \otimes T$. This will be a direct consequence of the fact that $\nabla \rP = 0$ and the following lemma:
\begin{lemm} \label{grass:lemm}
If $\mu = \mu_U \otimes \mu_{\mc{L}}$ or $\mu = \nu_U \otimes \nu_{\mc{L}}$ for $\mu_{\mc{L}}$, $\mu_U$ symmetric forms and $\nu_{\mc{L}}$, $\nu_U$ alternating forms, then conjugation by $g$ preserves $\mc{A}$.
\end{lemm}
\begin{lproof}
$\mu_U \otimes \mu_{\mc{L}} = (\mu_U \otimes Id_{\mc{L}})(Id_U \otimes \mu_{\mc{L}})$. Those elements commute, and conjugation by each one of them evidently preserves $\mc{A}$. Same result for $g = \nu_U \otimes \nu_{\mc{L}}$.
\end{lproof}
And the metric $\mu$ on $M$ is evidently of the first form.

The existence of quaternionic Kh\"aler manifolds is well known; the existence of quaternionic Kh\"aler manifolds with
\be
R + \partial \rP \neq 0,
\ee
comes from the fact that $\mf{sp}(p,q) \cdot \mf{sp}(1)$ on its standard representation is not of Ricci-type (\cite{spein} \cite{meric} and \cite{holclass}), so not determine only by its Ricci-tensor. In these cases, $\onab$ has non-vanishing curvature, hence is non-flat.
\end{proof}

The same results hold for almost Grassmannian manifolds, for the other real and complex forms of $\mf{g}$, $\mf{sl}(4m,\mbb{R})$ and $\mf{sl}(2m,\mbb{C})$. In those cases, we need to use (Einstein) metrics $g$ with holonomy algebra in $\mf{sp}(4m - 2, \mbb{R}) \cdot \mf{sl}(2,\mbb{R})$ and $\mf{sp}(2m - 2, \mbb{C}) \cdot \mf{sl}(2,\mbb{C})$ (see \cite{holclass}), and apply the results of Lemma \ref{grass:lemm} for $\mu = \nu_1 \otimes \nu_2$.

\begin{prop}
There exists Einstein involutions for path geometries.
\end{prop}
\begin{proof}

Let $M$ be a manifold. Then if $\nabla$ is the Levi-Civita connection of a positive-definite metric Einstein metric on $M$ with positive Einstein coefficient, then the projective Cartan connection generated by $\nabla$ has an Einstein involution $\sigma$ on it, generated by a positive-definite metric $g$ on $\mc{T}$ (see \cite{mepro1}). Since $g$ is positive-definite, $\sigma$ is a Cartan involution, meaning that the bilinear form
\be
B(\sigma - , - )
\ee
is of definite signature, for $B$ the Killing form on $\mc{A}$.

Paper \cite{CoresSpace} demonstrated that regular normal path geometries may be constructed from projective geometries. This gives a manifold map $N \to M$, with the path geometry Cartan connection on $N$ projecting down to the projective Cartan connection on $M$. Consequently, $\sigma$ lifts to a preserved involution on $\mc{A}^N$. Since $B(\sigma - , - )$ remains of positive signature, $B(\sigma(v),v) > 0$ for all local never zero sections $v$ of $TN^*$. Since $\mc{A}^N_{(0)} = (TN^*)^{\perp}$ via $B$, the condition
\be
\sigma(TN^*) \cap \mc{A}^N_{(0)} = 0
\ee
is automatic. Hence $\sigma$ is an Einstein involution for the path geometry.
\end{proof}

Similarly, we can extend Cartan involutions to correspondence spaces whenever they exist.

Putting all these results together, we have the theorem:
\begin{theo}
Examples of non-flat normal Tractor connections with Einstein involutions exist for conformal, projective, CR, contact-projective, path, almost quaternionic and degree two almost Grassmannian geometries.

Any Cartan involution on a geometry is automatically an Einstein involution and generates an Einstein involution on any of its correspondence spaces.
\end{theo}

General existence issues depend upon the existence of metrics $g$, complex structures $J$ or subbundles $K \subset \T$ satisfying the required properties and preserved by $\onab$. The general construction of such invariant structures is highly non-trivial. One way of constructing these may be to look at the formalism of BGG sequences (\cite{BGG} and \cite{BGG2}), which may allow us to build some examples of Einstein involutions. Unfortunately, though any Einstein involution will show up in the BGG sequence as the solution to an invariant differential operator, the converse is not true -- simply solving that invariant differential operator will not automatically produce an Einstein involution. However, these considerations are beyond the scope of this paper.

\subsubsection{Non-normal examples}

If we drop the normality condition, existence is trivial. For example, we may pick any splitting of $\mc{A} = T \oplus \mc{A}_0 \oplus T^*$ compatible with the projections of $\mc{A}$, and define $\onab$ as
\be
\nabla_X = X + \nabla_X + \rP(X).
\ee
We then call this splitting the one determined by the Weyl structure $\nabla$, and all the other Weyl structures and the splittings they define are determined by the action of $T^*$ on $\mc{A}$ (see \cite{TCPG}). Then if we've simply picked a $\rP$ with the required properties, and a connection $\nabla$ on $\mc{G}_0$ such that $\nabla \rP = 0$, we have generated a Tractor connection with a preserved Einstein involution. This construction has no real ties with the underlying geometry, but may be usefull for some existence results.

\section{Cone construction and the Einstein condition}
\subsection{The cone construction}

The cone construction is an attempt to generalise the projective cone construction \cite{mepro1} and the Einstein cone construction in confromal geometry to other settings. The idea is ultimately the same: calculate the holonomy of $\onab$ by replacing it with an equivalent affine connection $\wnab$, with minimal torsion and hence (hopefully) a holonomy group that is easier to calculate.
\begin{defi}[Cone construction] \label{cone:defi}
A cone construction for a parabolic manifold $(M^n, \mc{P}, \onab)$, and a Tractor bundle $\mc{V}$ is a manifold $N = \mc{C}(M)$ with an affine $\wnab$ and a submersion $\pi: N \to M$. These must obey the following conditions:
\begin{enumerate}
\item There is an action of the abelian group $\mbb{R}^n$ on $N$, such that the orbits of $\mbb{R}^n$ are exactly the fibres of $\pi$.
\item There is a subbundle ${\wt{\mc{V}}}$ of $TN$, preserved by $\wnab$, which has contains the vector fields generated by the $\mbb{R}^n$ action.
\item $\wnab$ is invariant under the action of $\mbb{R}^n$, hence so is ${\wt{\mc{V}}}$. This means that ${\wt{\mc{V}}}$ descends to a bundle on $M$, and $\wnab$ to a connection on that bundle.
\item ${\wt{\mc{V}}}/\mbb{R}^n = \mc{V}$ and there is an isomorphism $\wnab \cong \onab$ on $\mc{V}$.
\item The holonomies of $\onab$ on $\mc{V}$ and $\wnab$ on the ${\wt{\mc{V}}}$ are the same.
\item If $Q$ is a vertical vector field on $N$, and $Tor$ is the torsion of $\wnab$ on ${\wt{\mc{V}}}$, then $Tor(Q,-) = 0$. Consequently for vector fields $X$ and $Y$ on $M$, $Tor(X,Y)$ is well defined independently of the lifts of $X$ and $Y$.
\end{enumerate}
\end{defi}

\begin{defi}[Cone decomposition]
Let $\mc{V}$ be any Tractor bundle on which $\mc{A}$ acts faithfully. If it admits a decomposition as
\be
\mc{V} = \mc{V}_1 \oplus \mc{V}_2 \oplus \ldots \oplus \mc{V}_m,
\ee
such that each $\mc{V}_j$ is preserved by $\onab$ and admits a cone construction, we say that $(M,\mc{P}, \onab)$ admits a \emph{cone decomposition}. The point of the cone decomposition is that it allows us to use the cone construction in cases where $\mc{V}_1$ is preserved by $\onab$ but is not a Tractor bundle (which is the case for most reduced holonomy examples we have been looking at).
\end{defi}

These are the criteria for general cone construction, see for instance the projective cone \cite{mepro1}, which \emph{always} exists. Also relevant is the conformal double cone construction in conformal geometry \cite{dcone}. The double cone is instructive, as it is a cone only for a preferred connection with vanishing Cotton-York tensor. There are reasons to hope that a properly Einstein connection in the sense of Theorem \ref{main:theo} will make the existence of a cone more likely -- since such a connection will trivially satisfy any Cotton-York-like conditions, suppresses the difference between $T$ and its dual and may result in a decomposition of certain Tractor bundles into simpler pieces (thus making a cone decomposition more likely). Moreover, the preferred connection generated by the Einstein involution preserves a volume form on $T$, thus suppressing the distinction between density bundles for $T$ and sections of $\mbb{R} \times M$. The examples of the next section demonstrate that this is indeed the case.

\begin{theo}\label{cone:theo}
Given a parabolic geometry $(M, \mc{P}, \onab)$ with a preserved Einstein involution $\sigma$, let $\nabla$ be the Einstein connection generated by $\sigma$ as above. Let $\mc{B} \subset \A_0 = \A_0 \cap F_+ = \A_0 \cap \mf{so}(\rP)$ (note that this means that the holonomy subalgebra bundle of $\nabla$ is contained in $\mc{B}$, as $\nabla \rP = 0$). Assume that one of the Tractor bundles $\mc{V}$ of splits into $\mc{B}$-preserved components, as
\be
H \oplus L,
\ee
such that $H \subset T$, a non-degenerate subbundle for the metric $\rP$. Further assume that $L = \mbb{R} \times M$, that $H = \oplus_k H_k$, for distinct $\mc{B}$-irreducible bundles $H_k$, and that the action of $T^* \oplus T$ maps $L$ surjectively to $H$.

Then there exists a cone construction for $\mc{C}(M)$, the total space of the bundle $L$.
\end{theo}
\begin{proof}
Let $\pi$ be the projection $\mc{C}(M) \to M$; the $\mbb{R}$ action is evident. Define let $i$ as the embedding $T \subset T \oplus T^*$, sending the section $X$ of $T$ to
\be
X + \rP(X).
\ee
Now $T$ and $T^*$ are isomorphic under the action of $\mc{B}$, and $\rP$ sends $\mc{B}$-irreducible components to $\mc{B}$-irreducible components. So for each $H_k \subset \mc{V}$ there exists an $H_k \subset T$ such that $i(H_k)$ maps $L$ to $H_k$. Now the inclusion $H_k \subset T$ is only defined up to scale (as we have implicitly used a preserved volume form to get this inclusion). Now set $\eta$ as the constant section $1 \times M$ of $L$, and scale $H_k$ by the required amount to ensure that $i(X_k) \cdot \eta = X_k$ for $X_k$ a section of $H_k$.

Now $TL$ splits into vertical and horizontal vectors, $TL = W \oplus TM$, with $W \cong \mbb{R} \times L$. We then define $\wt{\mc{V}}$ as $H \oplus W$. Since $\mc{V} \cong H \oplus L$, there is a one-to-one correspondence between $\mbb{R}$-invariant sections of $\wt{\mc{V}}$ and sections of ${\mc{V}}$.

Under this identification, the connection $\wnab$ is:
\be
\wnab_X Z = 0 \ \ \textrm{and} \ \ \wnab_X \tau = \onab_{X_H} \tau, \ \ \textrm{and} \ \ \wnab_{\eta} X = X_H \ \ \textrm{and} \ \ \wnab_{\eta} \eta = \eta,
\ee
for $\mbb{R}$-invariant sections $X$ of $TM$, $Z$ of $H^{\perp}$, $\tau$ of $\mc{V}$. Here $X_H$ is defined as the image of $X$ under the orthogonal projection $TM \to H$, and we have extended $\eta$ into a section of $W$ by $\mbb{R}$-invariance. Since $\nabla$ acts trivially on $L$, we can see that
\be
\wnab_{X} \eta = \nabla_X \eta + i(X_H)\cdot \eta = X_H = \wnab_{\eta} X.
\ee
Since $[X,\eta] = -\mc{L}_{\eta} = 0$ as $X$ is $\mbb{R}$-invariant, this demonstrates that $\wnab$ has no torsion terms on vertical vectors. Also, $\wnab$ is evidently $\mbb{R}$-invariant, and preserves $\wt{\mc{V}}$, so we have all the properties of Definition \ref{cone:defi}, apart from the holonomy condition. We thus need to show that no extra holonomy is introduced. To do so, it suffices to demonstrate that:

\begin{lemm}
If $\tau$ is any section of $\mc{V}$, then it may be locally extended to a $\tau'$ of $\wt{\mc{V}}$ such that $\wnab_{\eta} \tau' = 0$, and, if $X$ is a section of $TM$ extended into $\mc{C}(M)$ $\mbb{R}$-invariantly, then
\be
\wnab_{\eta} \wnab_X \tau' = 0.
\ee
\end{lemm}
\begin{lproof}
Define a local coordinate $x$ on $\mc{C}(M)$ such that $\eta \cdot x = 1$, $X\cdot x = 0$ for $\mbb{R}$-invariant sections of $TM$, and $x = 0$ on $M\times 0$. Then set $\tau' = e^{-x} \tau$. This gives $\wnab_{\eta} \tau' = \tau' - \tau' = 0$. Moveover, since $X\cdot x = 0$, $\wnab_X e^{-x} \tau = e^{-x} \wnab_X \tau$, proving the lemma.
\end{lproof}

The preceding lemma shows that parallel transport along the $\eta$ direction will scale all vector fields by the same amount, but that this scaling commutes with parallel transported along directions in $M$. Since any loop $\psi$ starts and ends at the same point, this scaling will be cancel upon parallel transport along a $\psi$, and the result will be the same as if the vector had been parallel transported along the projected loop $\pi(\psi)$ in $M \times 0 \cong M$. Thus the holonomies of $\wnab$ on $\wt{\mc{V}}$ and $\onab$ on $\mc{V}$ are equal.

\end{proof}

\begin{rem}
There exists similar theorems when $L$ is not required to be a line-bundle, but they require extra technical conditions and are not needed for the known examples of cone constructions.
\end{rem}

\subsection{Examples of cones with the Einstein condition}
Here we present four examples of cone constructions or cone decompositions allowed by the Einstein condition.

\setcounter{exa}{0}
\begin{exa}
In conformal geometry of signature $(p,q)$, an Einstein preferred connection $\nabla$ implies a preserved section $\tau$ of the Tractor bundle. This means that the holonomy of $\onab$ acts faithfully on $\tau^{\perp}$. But if $\nabla$ is not Ricci-flat,
\be
\tau^{\perp} \cong T \oplus \mbb{R},
\ee
and it is easy to see that $T\oplus T^*$ maps $\mbb{R}$ surjectively to $T$. Moreover, here $F_+$ is an algebra bundle modelled on $\mf{so}(p,q+1)$ or $\mf{so}(p+1,q)$, and $\mc{B}$ is modelled on $\mf{so}(p,q)$. So all the hypotheses of Theorem \ref{cone:theo} are fulfilled, allowing a cone construction. Similarly, if a non-degenerate subbundle $K$ of $\mc{T}$ is preserved by $\onab$, this allows a cone decomposition \cite{mecon}.
\end{exa}

\begin{exa}
Let $(M,\onab)$ be a manifold with a free $m$-distribution (\cite{bryskew} and \cite{meskewnew}). In terms of the representation of parabolic geometries as Dynkin diagrams with crossed nodes \cite{two}, this geometry is given as $\skewnil$. It is alternately defined by a maximally non-integrable distribution $H$ of rank $m$. Then the standard Tractor bundle splits as $\mc{T} = H \oplus \mbb{R} \oplus H^*$, with a natural metric $h$ on it. Paper \cite{meskew} demonstrates that any Einstein involution must correspond to two preserved subbundles of $\mc{T}$:
\be
K \cong H \oplus \mbb{R} \ \ \textrm{and} \ \ K^{\perp} \cong H,
\ee
corresponding to $F_+$ being an algebra bundle modelled on $\mf{so}(p,q) \oplus \mf{so}(q,p+1)$ where $p+q = m$, and $\mc{B}$ being modelled on $\mf{so}(p,q)$. So each of these preserved pieces allow a cone construction: the $H$ piece has a trivial cone construction with $N = M$, and $\mbb{R}$ is mapped surjectively to $H$ as it is of non-zero norm and $\mf{so}(q,p+1)$ is surjective from any element onto its orthogonal complement. This thus gives $M$ a cone decomposition.
\end{exa}

\begin{exa}
Path geometries are given, in terms of the representation of parabolic geometries as Dynkin diagrams with crossed nodes \cite{two}, as $\pathgeo$. They have have $\mf{g} = \mf{sl}(n+2)$, $\mf{p} = (\mbb{R} \oplus \mf{gl}(n)) \rtimes (\mbb{R} \oplus \mbb{R}^n) \rtimes \mbb{R}^n$. They are characterised by subbundles $H_{-1} \oplus H_{-1}' \subset T$, for $H_{-1}'$ a line bundle, and $H_{-2} = H_{-1} \otimes H_{-1}' = T / (H_{-1} \oplus H_{-1}')$. A choice of Weyl structure $\nabla$ gives a splitting
\be
T = H_{-1} \oplus H_{-1}' \oplus H_{-2}.
\ee
If $\nabla$ preserves a volume form, the Tractor bundle $\mc{T}$ splits as:
\be
\mc{T} = H_{-1} \oplus \mbb{R} \oplus H_{-1}'.
\ee
Then algebraic considerations imply that $T \oplus T^*$ must map $\mbb{R}$ surjectively onto the other components of $\T$ and $\mc{B}$ is modelled on $\mf{so}(p,q)$ with $p+q = n$, giving the requirements for a cone construction.
\end{exa}

\begin{exa}
If we are in the case of the normal, torsion-free, degree-two almost Grassmannian Einstein geometry of subsection \ref{invo:exam}, there is a cone construction identical to that described in \cite{quacone}. The construction is similar to that of Theorem \ref{cone:theo}, except we are using quaternionic density bundles in this case.
\end{exa}

\section{Future research}

To extend these results, existence or non-existence proofs are needed for Einstein involutions in all non-flat normal Cartan connections. Analysing the exceptional cases would be interesting as well.

However, other possible future research is to look at those cases where the crucial property $\sigma(T^*) \cap \mc{A}_{(0)}$ fails or where $\sigma$ is replaced by a degenerate endomorphism. The author's paper on the geometry of free $m$-distributions \cite{meskew} suggests that in this case we will get a weakening of the uniqueness condition for the Einstein connection $\nabla$, as well as a weakening of the $\nabla \rP = 0$ condition, in that $\nabla \rP$ will only be zero when restricted to a certain subbundle of $\otimes^3 T^*$. If $\sigma$ is degenerate, then $\rP$ will often be degenerate as well.

\bibliographystyle{ieeetr}
\bibliography{Einstein}

\textbf{Stuart Armstrong}

Fakult\"at f\"ur Mathematik

Universit\"at Wien

Nordbergstr. 15, 1090 Wien

Austria

\vspace{5mm}

Lobenhauerngasse 8/10

1170 Wien

Austria

\vspace{5mm}

stuart.armstrong@stx.oxon.org

\end{document}